\documentstyle[12pt,intlim]{amsart}

\def\txt#1{\quad\text{#1}\quad} 

\def\C{\Bbb C} 
\def\R{\Bbb R} 
\newmathalphabet*{\eusm}{eus}{m}{n}
\newcommand{\Ker}{\operatorname{Ker}}
\renewcommand{\Im}{\operatorname{Im}}

\theoremstyle{plain} 
\newtheorem{proposition}{Proposition}
\newtheorem{theorem}{Theorem}		 \newtheorem{example}{Example}

\theoremstyle{definition}
\newtheorem{definition}{Definition}	

\theoremstyle{remark}
\newtheorem{remark}{Remark}		
\newtheorem{problem}{Problem}

\title{Real and Complex Operator Ideals} 
\author[J.~Wenzel]{J\"org Wenzel} 
\thanks{This article originated from the author's masters thesis at
  the University of Jena written under the supervision of A.~Pietsch} 
\address{Mathematical Institute\\FSU Jena\\07740 Jena\\Germany} 
\email{\tt wenzel@@minet.uni-jena.de}

\begin{document}
\maketitle

\begin{abstract}
  The powerful concept of an operator ideal on the class of all Banach
  spaces makes sense in the real and in the complex case. In both
  settings we may, for example, consider compact, nuclear, or
  $2$--summing operators, where the definitions are adapted to each
  other in a natural way. This paper deals with the question whether or
  not that fact is based on a general philosophy. Does there exists a
  one--to--one correspondence between ``real properties'' and ``complex
  properties'' defining an operator ideal? In other words, does there
  exist for every real operator ideal a uniquely determined
  corresponding complex ideal and vice versa?

  Unfortunately, we are not abel to give a final answer. Nevertheless,
  some preliminary results are obtained. In particular, we construct for
  every real operator ideal a corresponding complex operator ideal and
  for every complex operator ideal a corresponding real one. However, we
  conjecture that there exists a complex operator ideal which can not be
  obtained from a real one by this construction.

  The following approach is based on the observation that every
  complex Banach space can be viewed as a real Banach space with an
  isometry acting on it like the scalar multiplication by the imaginary
  unit $i$.
\end{abstract}


\section{Preliminaries}
\label{introduction}

Let $X$ always denote a Banach space over the field of real numbers.
The letters $A$, $B$, $T$ and $S$ refer to linear operators between
real Banach spaces. The identity map of $X$ is denoted by $I_X$.

For $x\in X$ the canonical injections
\[ X\to X\oplus X:x\mapsto(x,o)\txt{and}X\to X\oplus X:x\mapsto(o,x)
\] are denoted by $J_1^X$ and $J_2^X$, respectively.

For $x_1,x_2\in X$ the canonical surjections
\[ X\oplus X\to X:(x_1,x_2)\mapsto x_1\txt{and}X\oplus X\to
X:(x_1,x_2)\mapsto x_2
\] are denoted by $Q_1^X$ and $Q_2^X$, respectively.

Let $\boldsymbol{\frak L}$ always denote the ideal of all (real or
complex) bounded linear operators.

For the theory of operator ideals we refer to the monographs of
Pietsch, \cite{PIE:1} and \cite{PIE:2}.

\begin{definition}
  The Banach space $X\oplus X$ becomes a complex Banach space under
  the operations
  \begin{align*} 
    (x_1,x_2)+(y_1,y_2) &:= (x_1+y_1,x_2+y_2)\\
    (\alpha+i\beta)(x_1,x_2)&:= (\alpha x_1-\beta x_2,\beta x_1+\alpha
    x_2),
  \end{align*} 
  and the norm
  \[ \|(x_1,x_2)\|:=\left(\frac1{2\pi}\int_{-\pi}^{+\pi}\|x_1\cos\phi+x_2\sin
  \phi\|^2\,d\phi\right)^{1/2},
  \] where $x_1,x_2,y_1,y_2\in X$.

  The obtained complex Banach space is called the {\em
  complexification} of $X$. 
\end{definition}


\section{Banach spaces with an $i$--operator}

Note that the scalar multiplication by the imaginary unit $i$ in a
complex Banach space yields an isometric operator on this Banach
space. This observation leads to the following definition.

\begin{definition}
  A {\em Banach space with an $i$--operator} is a pair $[X,A]$
  consisting of a real Banach space $X$ and a linear operator $A:X\to X$
  such that $A^2=-I_X$ and
  \begin{equation}
    \|\alpha x+\beta Ax\|=\|x\|\quad\mbox{for}\
    |\alpha|^2+|\beta|^2=1.
    \label{conj-norm}
  \end{equation}
  We refer to the operator $A$ as an {\em $i$--operator\/} on the
  Banach space $X$.
\end{definition}

If $[X,A]$ is a Banach space with an $i$--operator, then we obviously
get
\[ (-A)^2=-I_X
\] and
\[ \|\alpha x-\beta A x\|=\|x\|\quad\mbox{for }|\alpha|^2+|\beta|^2=1,
\] and hence $[X,-A]$ is a Banach space with an $i$--operator, too.

\begin{definition}
  We denote by $\overline{[X,A]}:=[X,-A]$ the {\em complex conjugate
  Banach space} of $[X,A]$.
\end{definition}

A linear operator between two complex Banach spaces obviously commutes
with the scalar multiplication by $i$. The next definition is due to
this fact.

\begin{definition}
  Given two Banach spaces with $i$--operators $[X,A]$ and $[Y,B]$, we
  say that an operator $T:X\to Y$ respects $A$ and $B$ if
  \[ TA=BT.
  \] In this case, we write $[T,A,B]$ for the induced operator from
  $[X,A]$ to $[Y,B]$.
\end{definition}

If $T$ respects the $i$--operators $A$ and $B$, then we obviously get
\[ T(-A)=(-B)T
\] and hence $T$ respects the $i$--operators $-A$ and $-B$, too.

\begin{definition}
  We denote by
  \[ \overline{[T,A,B]}:=[T,-A,-B]:[X,-A]\to[Y,-B]
  \] the {\em complex conjugate operator} of $[T,A,B]$.
\end{definition}

The class of all Banach spaces with an $i$--operator $[X,A]$ together
with the morphisms $[T,A,B]$ forms a category. Two Banach spaces with
an $i$--operator $[X,A]$ and $[Y,B]$ are said to be isomorphic, if
they are isomorphic in the sense of this category. We then write
$[X,A]\simeq[Y,B]$.

The mappings
\[ [X,A]\to\overline{[X,A]}\quad\mbox{and}\quad[T,A,B]\to\overline{[T,A,B]}
\] yield a covariant functor on this category.


\section{Examples of Banach spaces with an $i$--operator}

\subsection{Complex Banach spaces}

Of course the concept of a Banach space with an $i$--operator exactly
simulates that of complex Banach spaces. This means if $Y$ is a
complex Banach space, then $[Y,A]$ is a Banach space with an
$i$--operator, provided $Y$ is considered as a real space and $A$ is
given by
\[ Ay:=i\cdot y\txt{for all $y\in Y$.}
\] In this way every complex Banach space can be considered as a
Banach space with an $i$--operator.

Also, if $Y$ and $Z$ are complex Banach spaces and $T:Y\to Z$ is a
complex operator then $T$ respects $A$ and $B$ defined as above on $Y$
and $Z$, respectively.

Conversely, given a Banach space with an $i$--operator $[X,A]$ we can
introduce a complex scalar multiplication in $X$ by
\[ i\cdot x:=Ax
\] or more generally
\[ (\alpha+i\beta)x:=\alpha x+\beta Ax\quad\mbox{for
  }\alpha,\beta\in\R.
\] In this way we obtain a complex Banach space.

Also, if $[T,A,B]$ is an operator respecting $A$ and $B$, then $T$
yields a com\-plex operator on the assigned complex spaces.

Thus the category of Banach spaces with an $i$--operator can be
identified with the category of complex Banach spaces.


\subsection{Complexification}\label{complexification}

For a real Banach space $X$, we can define a {\em natural
$i$--operator} $N_X$ on the space $X\oplus X$ by
\[ N_X(x_1,x_2):=(-x_2,x_1),\txt{for $x_1,x_2\in X$.}
\] It is easy to see that the space $[X\oplus X,N_X]$ equipped with
the norm
\[ \|(x,y)\|:=\left(\frac1{2\pi}\int_{-\pi}^{+\pi}\|x\cos\phi+
  y\sin\phi\|^2d\phi\right)^{1/2}
\] corresponds to the complexification of the space $X$ as defined in
Section \ref{introduction}.

Also, if $T:X\to Y$ is a real operator, then
\[ T\oplus T:(x_1,x_2)\mapsto(Tx_1,Tx_2)
\] respects the $i$--operators $N_X$ and $N_Y$. The corresponding
operator $[T\oplus T,N_X,N_Y]$ is called the {\em complexification} of
the operator $T$.

It is natural to ask whether or not every complex Banach space can be
obtained as the complexification of a real Banach space.  The
following proposition yields a characterization of complex Banach
spaces with that property. Example \ref{ex1} will then show that there
are complex Banach spaces which are not a complexification.

\begin{proposition}\label{prop1}
  The following conditions are equivalent:
  \begin{enumerate}
  \item The complex Banach space $[X,A]$ is a complexification of a
    real Banach space $Y$, i.e. there exists an isomorphism
    \[ [S,A,N_Y]:[X,A]\to[Y\oplus Y,N_Y]
    \] of the spaces $[X,A]$ and $[Y\oplus Y,N_Y]$. \label{1}
  \item There exists an automorphism $T$ of $X$ such that
    \[ TA=-AT\quad\mbox{and}\quad T^2=I_X,
    \] in particular the spaces $[X,A]$ and $\overline{[X,A]}$ are
    (complex) isomorphic. \label{2}
  \end{enumerate}
\end{proposition}

\begin{pf}
  First of all, note that the complex conjugation
  \[ C:Y\oplus Y\to Y\oplus Y:(y_1,y_2)\mapsto (y_1,-y_2)
  \] respects the $i$--operators $N_Y$ and $-N_Y$ and hence yields an
  isomorphism $[C,N_Y,-N_Y]$ of the complex conjugate spaces $[Y\oplus
  Y,N_Y]$ and $\overline{[Y\oplus Y,N_Y]}$.

  Moreover, for a given isomorphism $[S,A,N_Y]$ the complex conjugate
  map
  \[ [S,-A,-N_Y]=\overline{[S,A,N_Y]}:
    \overline{[X,A]}\to\overline{[Y\oplus Y,N_Y]}
  \] yields an isomorphism of the spaces $\overline{[X,A]}$ and
  $\overline{[Y\oplus Y,N_Y]}$.
  
  Putting together these isomorphisms, we get an isomorphism
  \[ [S^{-1}CS,A,-A]:[X,A]\to\overline{[X,A]}
  \] of the spaces $[X,A]$ and $\overline{[X,A]}$ which moreover
  satisfies
  \[ (S^{-1}CS)^2=I_X.
  \] This proves that \ref{1} implies \ref{2}.
  
  To prove the converse, we consider the subspace $Y\subseteq X$
  defined by
  \[ Y:=\{x+Tx:x\in X\}.
  \] Then the map
  \[ S:X\to Y\oplus Y:x\mapsto(Ax+TAx,x+Tx)
  \] respects the $i$--operators $A$ and $N_Y$. Moreover, by
  \eqref{conj-norm}, we have
  \begin{align*}
    \|Sx\| &= \left(\frac1{2\pi}\int_{-\pi}^{+\pi}
    \left\|(I+T)(Ax\cos\phi+x\sin\phi)\right\|^2d\phi\right)^{1/2}\\ &\leq
    \|I+T\|\left(\frac1{2\pi}\int_{-\pi}^{+\pi}
    \left\|Ax\cos\phi+x\sin\phi\right\|^2d\phi\right)^{1/2}\\ &=
    \|I+T\|\left(\frac1{2\pi}\int_{-\pi}^{+\pi}
    \|x\|^2d\phi\right)^{1/2}=\|I+T\|\cdot\|x\|
  \end{align*}
  and hence
  \[ \|S\|\leq\|I+T\|.
  \] By $T^2=I_X$ and $TA=-AT$, it can easily be verified that the
  inverse map of $S$ is given by
  \[ S^{-1}:Y\oplus Y\to X:(x_1+Tx_1,x_2+Tx_2)
    \mapsto\frac12\,\Big((x_2+Tx_2)-A(x_1+Tx_1)\Big).
  \] Consequently, $S$ defines an isomorphism of the spaces $[X,A]$
  and $[Y\oplus Y,N_Y]$.
\end{pf}


\subsection{Pathological examples}
\label{examples}

Note that on a finite dimensional space an $i$--operator can be
defined if and only if the space is of even dimension. However, there
are also examples of infinite dimensional Banach spaces not admitting
any $i$--operator.

\begin{example}[James \cite{jam51}]\label{ex2}
  The James space $J$ satisfies
  \[ \dim(J^{\ast\ast}/J)=1
  \] and hence cannot admit any $i$--operator.
\end{example}

\begin{remark}
  In {\rm \cite{sza86b}} Szarek gives an example of a superreflexive
  infinite dimensional Banach space that does not admit an
  $i$--operator.
\end{remark}

At first sight, one could think that the spaces $[X,A]$ and
$\overline{[X,A]}$ are quite similar. However, the following example
shows that they can differ in the worst imaginable way.

\begin{example}[Bourgain \cite{Bou86}]\label{ex1}
  There exists a Banach space with an $i$--operator $[X,A]$ such that
  $[X,A]$ and $\overline{[X,A]}$ are not isomorphic.
\end{example}

In view of the considerations in the next section the following
problem arises.

\begin{problem}\label{prob real squares}
  Does there exist a Banach space with an $i$--operator $[X,A]$ such
  that $X$ is not (real) isomorphic to a square of a real Banach space
  $Y$?
\end{problem}

By Proposition \ref{prop1} every complex Banach space $[X,A]$ which is
not isomorphic to $\overline{[X,A]}$ can not be (complex) isomorphic
to a square $Y\oplus Y$ equipped with the natural $i$--operator $N_Y$.
However, this problem only deals with real isomorphisms on real
spaces.


\subsection{Some relations}
The previous considerations show that there is a one--to--one
correspondence between the categories of complex Banach spaces and
real Banach spaces with an $i$--operator. Under this correspondence,
those complex Banach spaces that are a complexification of a real
space are exactly the Cartesian squares of real Banach spaces
$[X\oplus X,N_X]$ with the natural $i$--operator.

On the other hand, there is the forget functor acting from either the
category of complex Banach spaces or the category of real Banach
spaces with an $i$--operator into the category of real Banach spaces.
This functor simply forgets about the complex structure on a Banach
space. However, as the examples in the previous section show, this
functor is neither injective nor surjective.

Namely, the spaces $[X,A]$ and $\overline{[X,A]}$ of Example \ref{ex1}
are both mapped to the same real Banach space $X$, which shows that
the forget functor is not injective. The spaces of Example \ref{ex2}
show that the forget functor is not surjective.

Moreover, Proposition \ref{prop1} shows that the space $[X,A]$ of
Example \ref{ex1} is not a complexification of a real space.

\vspace{1cm}
\noindent
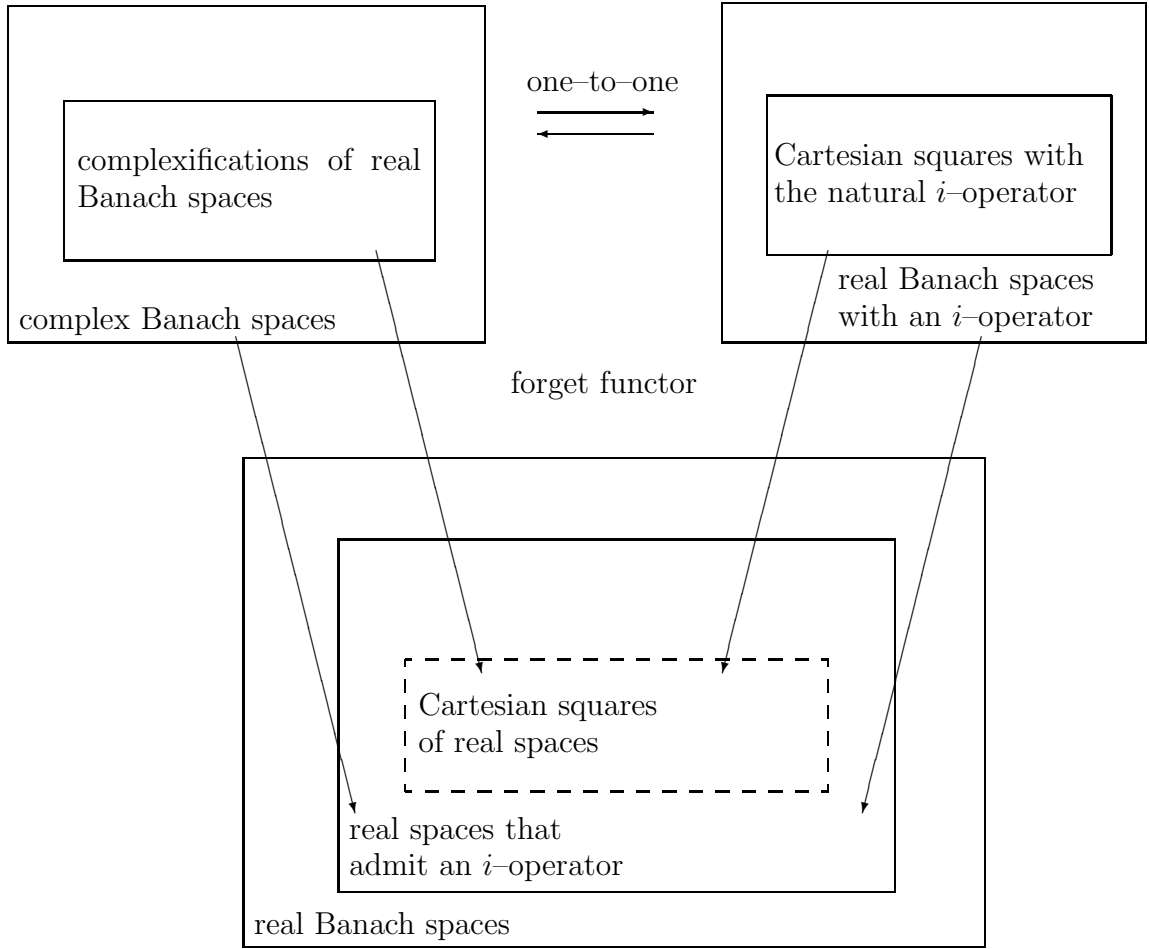
\begin{center}
  \unitlength=1.0pt
  \begin{picture}(430.00,384.00)(75.00,429.00)
    \put(214.00,713.00){\vector(1,-4){40.00}}
    \put(385.00,713.00){\vector(-1,-4){40.00}}
    \put(161.00,680.00){\vector(1,-4){45.00}}
    \put(443.00,680.00){\vector(-1,-4){45.00}}
    \put(275.00,765.00){\vector(1,0){44.00}}
    \put(319.00,757.00){\vector(-1,0){44.00}}
    \put(96.00,709.00){\framebox(140.00,60.00)[cc]{\parbox{130pt}
        {complexifications of real Banach spaces}}}
    \put(362.00,711.00){\framebox(130.00,60.00)[cc]{\parbox{125pt}
        {Cartesian squares with\\the natural $i$--operator}}}
    \put(75.00,678.00){\framebox(180.00,115.00)[bl]{\
        \rule[-5pt]{0pt}{30pt}\parbox[c]{170pt} {complex
          Banach\rule[-5pt]{0pt}{30pt} spaces}}}
    \put(345.00,678.00){\framebox(160.00,115.00)[bl]{\
        \rule[-13pt]{0pt}{30pt}\hspace*{1.4cm}\parbox[c]{150pt} {real Banach
          spaces \\with an $i$--operator}}}
    \put(225.00,508.00){\dashbox{5}(160.00,50.00)[cc]{\parbox{150pt}
        {Cartesian squares\\ of real spaces}}}
    \put(200.00,470.00){\framebox(210.00,120.00)[bl]{\
        \rule[-13pt]{0pt}{20pt}\parbox[c]{200pt} {real spaces that\\admit an
          $i$--operator}}} \put(164.00,449.00){\framebox(280.00,180.00)[bl]{\
        \rule[-5pt]{0pt}{20pt}\parbox[c]{270pt} {real Banach spaces}}}
    \put(271.00,777.00){\makebox(0,0)[lc]{one--to--one}}
    \put(265.00,662.00){\makebox(0,0)[lc]{forget functor}}
  \end{picture}
\end{center}

Problem \ref{prob real squares} deals with the question wether there
exists a real Banach space that admits an $i$--operator but is no real
square.


\section{Operator ideals}

In the following let $\boldsymbol{\frak R}$ and $\boldsymbol{\frak C}$
denote ideals of real or complex operators, respectively.

\begin{definition}
  We define the {\em complexification} $\boldsymbol{\frak R}_{\C}$ of
  the real operator ideal $\boldsymbol{\frak R}$ as the ideal consisting
  of all complex operators $[T,A,B]$ such that $T$ belongs to
  $\boldsymbol{\frak R}$, i.e.
  \[ \boldsymbol{\frak
    R}_{\C}([X,A],[Y,B]):=\left\{[T,A,B]\in\boldsymbol{\frak
    L}([X,A],[Y,B]):T\in\boldsymbol{\frak R}(X,Y)\right\}.
  \] We define the {\em real form} $\boldsymbol{\frak C}_{\R}$ of the
  complex operator ideal $\boldsymbol{\frak C}$ as the ideal consisting
  of all real operators $T$ such that its complexification belongs to
  $\boldsymbol{\frak C}$, i.e.
  \[ \boldsymbol{\frak C}_{\R}(X,Y):=\left\{T\in\boldsymbol{\frak
    L}(X,Y):[T\oplus T,N_X,N_Y]\in\boldsymbol{\frak C} ([X\oplus
    X,N_X],[Y\oplus Y,N_Y])\right\}.
  \]
\end{definition}

Easy computations show that for a given operator ideal
$\boldsymbol{\frak C}$ the class of all conjugate operators
$\overline{[T,A,B]}$ of operators $[T,A,B]\in\boldsymbol{\frak C}$ is
an operator ideal, too.

\begin{definition}
  For a complex operator ideal $\boldsymbol{\frak C}$, we denote by
  \[ \overline{\boldsymbol{\frak
      C}}:=\{\overline{[T,A,B]}:[T,A,B]\in\boldsymbol{\frak C}\}
  \] the {\em complex conjugate operator ideal} of $\boldsymbol{\frak
  C}$.
  
  We say that an operator ideal $\boldsymbol{\frak C}$ is self
  conjugate if $\overline{\boldsymbol{\frak C}}=\boldsymbol{\frak C}$.
\end{definition}

Of course one may ask whether or not the ideals $\boldsymbol{\frak R}$
and $\boldsymbol{\frak C}$ are uniquely determined by their complex or
real analogues $\boldsymbol{\frak R}_{\C}$ and $\boldsymbol{\frak
C}_{\R}$, respectively. This makes it necessary to examine the
combination of the procedures above, namely under which conditions on
$\boldsymbol{\frak R}$ and $\boldsymbol{\frak C}$ is it true that
\begin{equation}
  (\boldsymbol{\frak R}_{\C})_{\R}=\boldsymbol{\frak R}
  \label{real}
\end{equation}
and
\begin{equation}
  (\boldsymbol{\frak C}_{\R})_{\C}=\boldsymbol{\frak C}\,?
  \label{complex}
\end{equation}
Theorem \ref{theo reell} below will give a satisfactory answer if we
start from a real ideal $\boldsymbol{\frak R}$. In the complex case,
to fulfill \eqref{complex}, we must additionally require that the
ideal $\boldsymbol{\frak C}$ is self conjugate.


\section{The real case}

\begin{proposition}\label{real cart}
  A real operator ideal $\boldsymbol{\frak R}$ is uniquely determined
  by its restriction to Cartesian squares, i.e.
  \[ T\oplus T\in\boldsymbol{\frak R}(X\oplus X,Y\oplus Y)\iff
    T\in\boldsymbol{\frak R}(X,Y).
  \]
\end{proposition}

\begin{pf}
  Looking at the the diagrams
  \begin{center}
    \unitlength=1mm \linethickness{0.4pt}
    \begin{picture}(145.00,40.00)
      \put(14.00,30.00){\vector(1,0){42.00}}
      \put(17.00,10.00){\vector(1,0){36.00}}
      \put(10.00,26.00){\vector(0,-1){12.00}}
      \put(60.00,14.00){\vector(0,1){12.00}}
      \put(85.00,14.00){\vector(0,1){12.00}}
      \put(135.00,26.00){\vector(0,-1){12.00}}
      \put(89.00,30.00){\vector(1,0){42.00}}
      \put(92.00,10.00){\vector(1,0){36.00}}
      \put(35.00,33.00){\makebox(0,0)[cc]{$T$}}
      \put(35.00,13.00){\makebox(0,0)[cc]{$T\oplus T$}}
      \put(110.00,13.00){\makebox(0,0)[cc]{$T\oplus T$}}
      \put(110.00,33.00){\makebox(0,0)[cc]{$T$}}
      \put(6.00,20.00){\makebox(0,0)[cc]{$J_1^X$}}
      \put(64.00,20.00){\makebox(0,0)[cc]{$Q_1^Y$}}
      \put(81.00,20.00){\makebox(0,0)[cc]{$Q_1^X$}}
      \put(89.00,20.00){\makebox(0,0)[cc]{$Q_2^X$}}
      \put(131.00,20.00){\makebox(0,0)[cc]{$J_1^Y$}}
      \put(139.00,20.00){\makebox(0,0)[cc]{$J_2^Y$}}
      \put(135.00,30.00){\makebox(0,0)[cc]{$Y$}}
      \put(85.00,30.00){\makebox(0,0)[cc]{$X$}}
      \put(60.00,30.00){\makebox(0,0)[cc]{$Y$}}
      \put(10.00,30.00){\makebox(0,0)[cc]{$X$}}
      \put(10.00,10.00){\makebox(0,0)[cc]{$X\oplus X$}}
      \put(60.00,10.00){\makebox(0,0)[cc]{$Y\oplus Y$}}
      \put(85.00,10.00){\makebox(0,0)[cc]{$X\oplus X$}}
      \put(135.00,10.00){\makebox(0,0)[cc]{$Y\oplus Y$,}}
    \end{picture}
  \end{center}
  we easily see that
  \[ T=Q_1^Y(T\oplus T)J_1^X\txt{and}T\oplus
  T=J_1^YTQ_1^X+J_2^YTQ_2^X.
  \] The assertion now follows from the ideal properties of
  $\boldsymbol{\frak R}$.
\end{pf}

The following theorem means that real operator ideals are uniquely
determined by their complexification.

\begin{theorem}\label{theo reell}
  For any real operator ideal $\boldsymbol{\frak R}$ we have
  \[ (\boldsymbol{\frak R}_{\C})_{\R}=\boldsymbol{\frak R}.
  \]
\end{theorem}

\begin{pf}
  Let $T\in(\boldsymbol{\frak R}_{\C})_{\R}$. By the definition of the
  real form of an operator ideal this is the case if and only if
  $[T\oplus T,N_X,N_Y]\in\boldsymbol{\frak R}_{\C}$. By the definition
  of $\boldsymbol{\frak R}_{\C}$ this is equivalent to $T\oplus
  T\in\boldsymbol{\frak R}$. Now, using Proposition \ref{real cart}, the
  assertion follows.
\end{pf}


\section{The complex case}

Let $[X,A]$ be a Banach space with an $i$--operator. Then $A\oplus -A$
yields an $i$--operator on $X\oplus X$. The following proposition
shows the significance of this $i$--operator.

\begin{proposition}\label{squares}
  Let $[X,A]$ be a Banach space with an $i$--operator. In the category
  of Banach spaces with an $i$--operator the space $[X\oplus X,A\oplus
  -A]$ is isomorphic to the complexification of the real space $X$,
  \[ [X\oplus X,A\oplus-A]\simeq[X\oplus X,N_X].
  \]
\end{proposition}

\begin{pf}
  Let $x_1,x_2\in X$ and define $T:X\oplus X\to X\oplus X$ by
  \[ T(x_1,x_2):=(x_1+Ax_2,x_1-Ax_2).
  \] Then we have
  \begin{align*}
    TN_X(x_1,x_2) &= (Ax_1-x_2,-Ax_1-x_2),\\ (A\oplus-A)T(x_1,x_2) &=
    (Ax_1-x_2,-Ax_1-x_2).
  \end{align*} 
  Hence the operator $T$ respects the $i$--operators $N_X$ and
  $A\oplus-A$. It is evident that $T$ is a bounded bijection and
  consequently $[T,N_X,A\oplus-A]$ defines an isomorphism of the spaces
  $[X\oplus X,N_X]$ and $[X\oplus X,A\oplus -A]$.
\end{pf}

Note that, in order to make $N_X$ an $i$--operator on $X\oplus X$, we
must modify the norm a little bit, introducing e.g.
\[ \|(x,y)\|:=\left(\frac1{2\pi}\int_{-\pi}^{+\pi}
\|x\cos\phi+y\sin\phi\|^2d\phi\right)^{1/2}.
\] On the other hand, in order to make $A\oplus A$ an $i$--operator on
$X\oplus X$, we may use e.g.
\[ \|(x,y)\|:=\|x\|+\|y\|.
\] Nevertheless, the resulting spaces are still isomorphic.

In contrast to Proposition \ref{real cart}, in the complex case, we
can only prove the following result.

\begin{proposition}\label{complex cart}
  For a complex operator ideal $\boldsymbol{\frak C}$, we have
  \begin{align*}
    [T\oplus T,A\oplus -A,B\oplus -B]&\in\boldsymbol{\frak C}([X\oplus
    X,A\oplus -A], [Y\oplus Y,B\oplus -B])\,\Longrightarrow\\{}
    [T,A,B]&\in \boldsymbol{\frak C}([X,A],[Y,B]).
  \end{align*}
\end{proposition}

\begin{pf}
  First of all, note that for a Banach space with an $i$--operator
  $[X,A]$ the canonical injection $J_1^X:X\to X\oplus X$ respects the
  $i$--operators $A$ and $A\oplus -A$ whereas the canonical injection
  $J_2^X:X\to X\oplus X$ respects the $i$--operators $-A$ and $A\oplus
  -A$.

  Similarly, the canonical surjection $Q_1^X:X\oplus X\to X$ respects
  the $i$--operators $A\oplus -A$ and $A$ whereas the canonical
  surjection $Q_2^X:X\oplus X\to X$ respects the $i$--operators
  $A\oplus-A$ and $-A$.
  
  Hence we can form the composition of the operators $[J_1^X,A,A\oplus
  -A]$, $[T\oplus T,A\oplus-A,B\oplus-B]$ and $[Q_1^Y,B\oplus-B,B]$ and
  we have
  \[ [T,A,B]=[Q_1^Y(T\oplus T)J_1^X,A,B].
  \] Now it follows from the ideal properties of $\boldsymbol{\frak
  C}$ that $[T,A,B]\in\boldsymbol{\frak C}$.
\end{pf}

To prove the converse of Proposition \ref{complex cart} an additional
property of the operator ideal $\boldsymbol{\frak C}$ is required.

\begin{proposition}\label{complex cart 2}
  Let $\boldsymbol{\frak C}$ be a complex operator ideal. The
  following conditions are equivalent:
  \begin{enumerate}
  \item\label{cond1} The operator ideal $\boldsymbol{\frak C}$ is self
    conjugate.
  \item\label{cond2} The operator ideal $\boldsymbol{\frak C}$ is
    uniquely determined by its restriction to Cartesian squares, i.e.
    \begin{align*}
      [T\oplus T,A\oplus -A,B\oplus -B]&\in\boldsymbol{\frak
      C}([X\oplus X,A\oplus -A], [Y\oplus Y,B\oplus -B])\iff\\{}
      [T,A,B]&\in\boldsymbol{\frak C} ([X,A],[Y,B]).
    \end{align*}
  \end{enumerate}
\end{proposition}

\begin{pf}
  If $\boldsymbol{\frak C}$ is self conjugate then it follows from
  $[T,A,B]\in\boldsymbol{\frak C} $ that $[T,-A,-B]\in\boldsymbol{\frak
  C}$. Hence we may form the composition of the operators
  $[Q_1^X,A\oplus-A,A]$, $[T,A,B]$ and $[J_1^Y,B,B\oplus-B]$ as well as
  of the operators $[Q_2^X,A\oplus-A,-A]$, $[T,-A,-B]$ and
  $[J_2^Y,-B,B\oplus-B]$. The proof of condition \ref{cond2} now follows
  from the formula
  \[ [T\oplus T,A\oplus-A,B\oplus-B]=[J_1^YTQ_1^X,A\oplus-A,B\oplus-B]
    + [J_2^YTQ_2^X,A\oplus-A,B\oplus-B]
  \] and the ideal properties of $\boldsymbol{\frak C}$.

  If on the other hand condition \ref{cond2} is satisfied, then it
  follows from $[T,A,B]\in\boldsymbol{\frak C}$ that
  \[ [T\oplus T,A\oplus-A,B\oplus-B]\in\boldsymbol{\frak C}.
  \] Now we may form the composition of $[J_2^X,-A,A\oplus-A]$,
  $[T\oplus T,A\oplus-A,B\oplus-B]$ and $[Q_2^Y,B\oplus-B,-B]$. The
  identity
  \[ [T,-A,-B]=[Q_2^Y(T\oplus T)J_2^X,-A,-B]
  \] and the ideal properties of $\boldsymbol{\frak C}$ imply that
  $[T,-A,-B]\in\boldsymbol{\frak C}$ and hence $\boldsymbol{\frak C}$ is
  self conjugate.
\end{pf}

Now we can prove the complex analogue to Theorem \ref{theo reell}.

\begin{theorem}\label{theo complex}
  For any complex operator ideal $\boldsymbol{\frak C}$ we have
  \[ (\boldsymbol{\frak C}_{\R})_{\C}\subseteq\boldsymbol{\frak C}.
  \]
  
  In addition, the following conditions are equivalent:
  \begin{enumerate}
  \item The operator ideal $\boldsymbol{\frak C}$ is self conjugate.
  \item\label{cond4} The operator ideal $\boldsymbol{\frak C}$ is
    uniquely determined by its real form, i.e.
    \[ (\boldsymbol{\frak C}_{\R})_{\C}=\boldsymbol{\frak C}.
    \]
  \end{enumerate}
\end{theorem}

\begin{pf}
  Let $[T,A,B]\in(\boldsymbol{\frak C}_{\R})_{\C}$. By the definition
  of the complexification of an operator ideal this is the case if and
  only if $T\in\boldsymbol{\frak C}_{\R}$. By the definition of
  $\boldsymbol{\frak C}_{\R}$ this is equivalent to $[T\oplus
  T,N_X,N_Y]\in \boldsymbol{\frak C}$. Next, using Proposition
  \ref{squares}, this is the same as to say that $[T\oplus
  T,A\oplus-A,B\oplus -B]\in\boldsymbol{\frak C}$. Now, from Proposition
  \ref{complex cart} we see that $[T,A,B]\in\boldsymbol{\frak C}$
  follows. This proves the first part of the theorem.
  
  If the operator ideal $\boldsymbol{\frak C}$ is self conjugate, then
  by Proposition \ref{complex cart 2} it follows from
  $[T,A,B]\in\boldsymbol{\frak C}$ that $[T\oplus T,A\oplus -A,B\oplus
  -B]\in\boldsymbol{\frak C}$ and this implies that
  $[T,A,B]\in(\boldsymbol{\frak C}_{\R})_{\C}$.
  
  If on the other hand
  \[ (\boldsymbol{\frak C}_{\R})_{\C}=\boldsymbol{\frak C},
  \] then $[T,A,B]\in\boldsymbol{\frak C}$ implies that
  $[T,A,B]\in(\boldsymbol{\frak C}_{\R})_{\C}$ and hence $[T\oplus
  T,A\oplus -A,B\oplus -B]\in\boldsymbol{\frak C}$.  Consequently, it
  follows from Proposition \ref{complex cart 2} that $\boldsymbol{\frak
  C}$ is self conjugate.
\end{pf}


\section{Self conjugate operator ideals}

Theorem \ref{theo complex} gives a good criterion to check whether or
not a complex operator ideal is uniquely determined by its real form.

However, if one examines known operator ideals, it turns out, that all
of them are self conjugate. Hence the problem arises whether there
exist operator ideals which are not self conjugate. No such example
seems to be known.

\begin{problem}
  Does there exist a complex operator ideal $\boldsymbol{\frak C}$ and
  an operator $[T,A,B]\in\boldsymbol{\frak C}$ such that
  $\overline{[T,A,B]}=[T,-A,-B]\notin\boldsymbol{\frak C }$?
\end{problem}

Proposition \ref{example} will reduce this problem to the following
problem of the existence of a certain Banach space.

\begin{problem}\label{prob2}
  Does there exist a complex Banach space $[X,A]$ such that the
  following two conditions are satisfied
  \begin{eqnarray}\label{eqn1}
    [X,A]&\not\simeq&\overline{[X,A]}=[X,-A]\\{}
    \label{eqn5}\label{eqn2}
    [X,A]&\simeq&[X\oplus X,A\oplus A]?
  \end{eqnarray}
\end{problem}

The examples given in Section \ref{examples} give an idea of how to
construct such a space.

\begin{proposition}\label{example}
  Let $[X,A]$ be a Banach space with the properties \eqref{eqn1} and
  \eqref{eqn2}.  Then the operator ideal $\boldsymbol{\frak C}$ defined
  by
  \[ \boldsymbol{\frak C}:=\{[T,A,B]:\mbox{\rm$[T,A,B]$ admits a
    factorization over the space $[X,A]$}\}
  \] is not self conjugate.
\end{proposition}

\begin{pf}
  First of all, we see from condition \eqref{eqn2} that
  $\boldsymbol{\frak C}$ is indeed an operator ideal. Now, on the one
  hand, we surely have
  \[ [I_X,A,A]\in\boldsymbol{\frak C}.
  \] Suppose that
  $\overline{[I_X,A,A]}=[I_X,-A,-A]\in\boldsymbol{\frak C}$. Then there
  exists a factorization
  \[ [I_X,-A,-A]=[R,A,-A][S,-A,A],
  \] where $R,S:X\to X$. Let $Y:=\Im(S)=\Im(SR)$ and $Z:=\Ker(SR)$.
  Then $Y$ is a com\-ple\-mented subspace in $X$ and therefore we have
  \begin{eqnarray}\label{eqn3}
    [X,A]&\simeq&[Y\oplus Z,A\oplus A]\\{}
    \label{eqn7}[X,-A]&\simeq&[Y\oplus Z,-A\oplus -A]\\{}
    \label{eqn8}[Y,A]&\simeq&[X,-A]\\{}\label{eqn4}
    [Y,-A]&\simeq&[X,A].
  \end{eqnarray}
  Moreover, it follows from \eqref{eqn5} that
  \begin{equation}
    [X,-A]\simeq[X\oplus X,-A\oplus -A].
    \label{eqn6}
  \end{equation}
  The following is also known as Pe\l czy\'nski's decomposition
  method.  Using \eqref{eqn3} through \eqref{eqn4}, we get

  \noindent
  \begin{tabular*}{\textwidth}{@{\hspace{5cm}}
r@{\extracolsep{\tabcolsep}}c@{\extracolsep{\tabcolsep}}l@{\extracolsep{\fill}}
l@{\hspace{0pt}}}
    $[X,A]$&$\simeq$&$[Y\oplus Z,A\oplus A]$&by \eqref{eqn3}\\
    &$\simeq$&$[X\oplus Z,-A\oplus A]$&by \eqref{eqn8}\\
    &$\simeq$&$[X\oplus X\oplus Z,-A\oplus -A\oplus A]$&by \eqref{eqn6}\\
    &$\simeq$&$[X\oplus Y\oplus Z,-A\oplus A\oplus A]$&by \eqref{eqn8}\\
    &$\simeq$&$[X\oplus X,-A\oplus A]$&by \eqref{eqn3}\\
    &$\simeq$&$[Y\oplus Z\oplus X,-A\oplus -A\oplus A]$&by \eqref{eqn7}\\
    &$\simeq$&$[X\oplus Z\oplus X,A\oplus -A\oplus A]$&by \eqref{eqn4}\\
    &$\simeq$&$[X\oplus Z,A\oplus -A]$&by \eqref{eqn5}\\
    &$\simeq$&$[Y\oplus Z,-A\oplus -A]$&by \eqref{eqn4}\\
    &$\simeq$&$[X,-A]$&by \eqref{eqn7}.
  \end{tabular*} 

  This yields a contradiction to property \ref{eqn1} of the space
  $[X,A]$ and therefore completes the proof.
\end{pf}




\end{document}